\pdfoutput=1
\RequirePackage{ifpdf}
\ifpdf 
\documentclass[pdftex]{sigma}
\else
\documentclass{sigma}
\fi

\input{xy}
\xyoption{all}

\usepackage{bm}
\def\P #1{\partial_{#1}}
\def\ch #1{\text{\rm Char}\ #1}

\begin{document}

\allowdisplaybreaks

\renewcommand{\thefootnote}{$\star$}

\renewcommand{\PaperNumber}{024}

\FirstPageHeading

\ShortArticleName{Cauchy Problem and Darboux Integrable Wave Maps}

\ArticleName{Cauchy Problem for a~Darboux Integrable\\
Wave Map System and Equations of Lie Type\footnote{This paper is a~contribution to the Special Issue ``Symmetries of
Dif\/ferential Equations: Frames, Invariants and~Applications''.
The full collection is available at
\href{http://www.emis.de/journals/SIGMA/SDE2012.html}{http://www.emis.de/journals/SIGMA/SDE2012.html}}}

\Author{Peter J.~VASSILIOU}

\AuthorNameForHeading{P.J.~Vassiliou}

\Address{Program in Mathematics and Statistics, University of Canberra, 2601 Australia}
\Email{\href{mailto:peter.vassiliou@canberra.edu.au}{peter.vassiliou@canberra.edu.au}}

\ArticleDates{Received September 27, 2012, in f\/inal form March 12, 2013; Published online March 18, 2013}

\Abstract{The Cauchy problem for harmonic maps from Minkowski space with its standard f\/lat metric to a~certain
non-constant curvature Lorentzian 2-metric is studied.
The target manifold is distinguished by the fact that the Euler--Lagrange equation for the ener\-gy functional is Darboux
integrable.
The time evolution of the Cauchy data is reduced to an ordinary dif\/ferential equation of Lie type associated to
${\rm SL}(2)$ acting on a~manifold of dimension~4.
This is further reduced to the simplest Lie system: the Riccati equation.
Lie reduction permits explicit representation formulas for various initial value problems.
Additionally, a~concise (hyperbolic) Weierstrass-type representation formula is derived.
Finally, a~number of open problems are framed.
}


\Classification{53A35; 53A55; 58A15; 58A20; 58A30}

\renewcommand{\thefootnote}{\arabic{footnote}}
\setcounter{footnote}{0}

\begin{flushright}
\begin{minipage}{105mm}
\it This paper is dedicated to Peter Olver on the occasion of his 60th birthday in celebration of his
contributions to mathematics; especially his influential, diverse applications of Lie theory.
\end{minipage}
\end{flushright}

\section{Introduction}

Let $(M,g)$ and $(N,h)$ be Riemannian or pseudo-Riemannian manifolds and $\varphi:M\to N$ a~smooth map.
The {\it energy} of $\varphi$ over a~compact domain $\mathcal{D}\subseteq M$ is
\begin{gather*}
e(\varphi)=
\frac{1}{2}\int_\mathcal{D}g^{ij}(x)h_{\alpha\beta}(\varphi)\frac{\partial\varphi^\alpha}{\partial x_i}\frac{\partial\varphi^\beta}{\partial x_j}d\text{vol}_M.
\end{gather*}
The critical points of $e(\varphi)$ satisfy the partial dif\/ferential equation (PDE)
\begin{gather*}
\triangle\varphi^\gamma+g^{ij}\Lambda^\gamma_{\alpha\beta}\frac{\partial\varphi^\alpha}{\partial x_i}\frac{\partial\varphi^\beta}{\partial x_j}=0,
\end{gather*}
where \looseness=-1 $\triangle$ is the Laplacian on $M$ and $\Lambda^\gamma_{\alpha\beta}$ the Christof\/fel symbols on $N$.
A map $\varphi$ is said to be {\it harmonic} if it is critical for $e(\varphi)$.
Harmonic maps generalise harmonic functions and geodesics and have been under intense study since the pioneering work
of Eells and Sampson~\cite{EelsSampson}; see~\cite{BairdWood} and~\cite{Helein} and references therein for
comprehensive introductions to the f\/ield.
If $\text{domain}(\varphi)=M=\mathbb{R}$ then harmonic maps are geodesic f\/lows.
If $\text{codomain}(\varphi)=N=\mathbb{R}$ then harmonic maps are harmonic functions.
If $(M,g)$ is pseudo-Riemannian then harmonic maps are called {\it wave maps}.
In this paper our focus is on wave maps, specif\/ically, the case $(M,g)=(\mathbb{R}^{1,1},dx\,dy)$.
Any further reference to wave maps in this paper means the domain space~$M$ is Minkowski space~$\mathbb{R}^{1,1}$ with
its standard f\/lat metric.
To provide slightly more insight, a~``physical'' illustration of wave maps in this class can be given: the motion of
a~frictionless elastic string constrained to vibrate on Riemannian surface~$(N,h)$, such as a~sphere, is exactly
described by a~wave map into~$N$; see~\cite{tao}.

There is a~well-known geometric literature on wave maps that has developed over the last two decades, especially their
existence as solutions of completely integrable systems; see~\cite{Guest} for a~textbook account with many references.
There is a~closely related physics literature where the relevant systems are known as nonlinear sigma models;
see~\cite{Zak}.
The f\/irst person to treat the Cauchy problem for wave maps into Riemannian targets was Chao-Hao Gu~\cite{Gu}.
He established the fundamental result that for smooth initial data, wave maps into complete Riemannian metrics have
long-time existence.
Gu's work initiated many further investigations where regularity constraints on the initial data have been
signif\/icantly relaxed.
Furthermore some higher dimensional problems have been treated; see~\cite{Struwe}.

In this paper we initiate the study of the Cauchy problem for wave maps in the special case where the systems they
satisfy are {\it Darboux integrable}.
Our f\/irst main result proves that the solution of the Cauchy problem for such a~Darboux integrable nonlinear sigma
model can be quite explicitly expressed as the f\/low of a~special vector f\/ield~$\xi_{k_1,k_2}$ (see
Theorem~\ref{IVPmain}) which itself is a~curve in a~certain Lie algebra of vector f\/ields canonically and
intrinsically associated to the Darboux integrable nonlinear sigma model, namely, its {\it Vessiot algebra}.
In consequence of this, standard constructions which facilitate the resolution of {\it systems of Lie type} such as
{\it Lie reduction} become available to the solution of the Cauchy problem for such wave map systems.
For this reason we have included an appendix to this paper which gives a~brief summary of the main results on systems
of Lie type adapted to the applications we envisage.
In this paper we have decided to focus on just one interesting nonlinear sigma model in order to discuss the
relationship between Darboux integrable hyperbolic systems on the one hand and the resolution of the corresponding
Cauchy problem via dif\/ferential systems of Lie type and to do so as explicitly as possible.
However, it will be seen that the proof of Theorem~\ref{IVPmain} is easy to generalise to other Darboux integrable
systems.
Indeed the very recent work~\cite{AndersonFels} outlines a~general, intrinsic proof of the close relationship between
systems of Lie type and the Cauchy problem for a~wide class of Darboux integrable exterior dif\/ferential systems.
The second main result of the paper marshalls the general theory of Darboux integrable exterior dif\/ferential
systems~\cite{AFV}, and generalised Goursat normal form~\cite{Vassiliou2006a,Vassiliou2006b} to derive
a~hyperbolic Weierstrass-type representation (Theorem~\ref{Wrep}) for wave maps into the non-constant curvature
metric~\eqref{MainMaps}.

As the name implies the notion of Darboux integrability originated in the 19th century and was most signif\/icantly
developed by Goursat~\cite{GoursatBook}.
Classically, it was a~method for constructing the ``general solution" of second order PDE in one dependent and two
independent variables
\begin{gather*}
F(x,y,u,u_x,u_y,u_{xx},u_{xy},u_{yy})=0
\end{gather*}
that generalised the so called ``method of Monge".
It relies on the notion of {\it characteristics} and their {\it first integrals}.
We refer the reader to~\cite{GoursatBook,IveyLandsbergBook,Vassiliou2001,Vassiliou2000} for
further information on classical Darboux integrability.
There are also extensive studies of Darboux integrable systems relevant to the equation class under study in the
works~\cite{Zhiber} and~\cite{ZhiberSokolov}.

\looseness=-1
In this paper we use a~new geometric formulation of Darboux integrable exterior dif\/ferential systems~\cite{AFV}.
At the heart of this theory is the fundamental notion of a~{\it Vessiot group} which, together with systems of Lie type
are our main tools for the study of the Cauchy problem for wave maps.

The PDE that govern wave maps have the semilinear form
\begin{gather*}
\boldsymbol{u}_{xy}=\boldsymbol{f}(x,y,\boldsymbol{u},\boldsymbol{u}_x,\boldsymbol{u}_y),
\qquad
\boldsymbol{u},\boldsymbol{f}\in\mathbb{R}^n.
\end{gather*}
Each solution possesses a~{\it double foliation} of curves called {\it characteristics}.
Such PDE often model {\it wave-like phenomena} and projection of these curves into the independent variable space
describes the space-time history of the wave propagation.
The characteristics of $\boldsymbol{u}_{xy}=\boldsymbol{f}$ are the integral curves of a~pair of rank $n+1$
distributions
\begin{gather*}
H_1=\left\{D_x+D_y\boldsymbol{f}\cdot\P{\boldsymbol{u}_{yy}},  \, \P{\boldsymbol{u}_{xx}}\right\},
\qquad
H_2=\left\{D_y+D_x\boldsymbol{f}\cdot\P{\boldsymbol{u}_{xx}}, \, \P{\boldsymbol{u}_{yy}}\right\},
\end{gather*}
where
\begin{gather*}
D_x=\P x+\boldsymbol{u}_x\P{\boldsymbol{u}}+\boldsymbol{u}_{xx}\P{\boldsymbol{u}_x}+\boldsymbol{f}\P{\boldsymbol{u}_y},
\qquad
D_y=\P y+\boldsymbol{u}_y\P{\boldsymbol{u}}+\boldsymbol{f}\P{\boldsymbol{u}_x}+\boldsymbol{u}_{yy}\P{\boldsymbol{u}_y},
\end{gather*}
are the {\it total differential operators} along solutions of the PDE.
Note that if $\boldsymbol{\theta}$ is the standard Cartan codistribution for $\boldsymbol{u}_{xy}=\boldsymbol{f}$ then
$H_1\oplus H_2=\text{ann}\,\boldsymbol{\theta}$.
Distributions $H_i$ are well-def\/ined with canonical structure.
\begin{definition}
If $\Delta$ is a~distribution on manifold $M$ then a~function $f:M\to\mathbb{R}$ is said to be a~{\it first integral}
of $\Delta$ if $Xf=0$ for all $X\in\Delta$.
\end{definition}

There is a~geometric def\/inition of Darboux integrable exterior dif\/ferential system~\cite{AFV}.
For wave map equations it reduces to
\begin{definition}
A semilinear system $\boldsymbol{u}_{xy}=\boldsymbol{f}$ with $\boldsymbol{u},\boldsymbol{f}\in\mathbb{R}^n$ is Darboux integrable
at a~given order if each of its characteristic systems $H_i$ has at least $n+1$ independent f\/irst integrals at that
order.
\end{definition}

\section{The Cauchy problem}\label{section2}

We consider wave maps
\begin{gather}
\label{MainMaps}
\boldsymbol{u}: \ \big(\mathbb{R}^{1,1},dx dy\big)\to\left(\mathbb{R}^2,\frac{du_1^2-du_2^2}{1+e^{-u_1}}\right).
\end{gather}
R.~Ream~\cite{Ream} studied the PDE for wave maps into nonzero curvature surface metrics that are Darboux integrable on
the 2-~and~3-jets and proved a~theorem that any such metric is (real) equivalent to one or other of the metrics
\begin{gather*}
\rho_\pm:=\frac{du_1^2+du_2^2}{1\pm e^{u_1}}.
\end{gather*}
Here \looseness=-1 we consider a~semi-Riemannian version of a~Ream metric and study the corresponding Cauchy problem.
We show how the solution of the Cauchy problem for wave maps~\eqref{MainMaps} can be expressed as an ordinary
dif\/ferential equation of {\it Lie type}.
Indeed we prove that the solution of the Cauchy problem for wave maps is naturally equivalent to an initial value
problem of a~Lie system for a~local action of ${\rm SL}(2)$ on a~manifold that is locally dif\/feomorphic to $\mathbb{R}^4$.
This is further reduced to an initial value problem for a~{\it single Riccati equation} together with {\it
a~quadrature}.

\looseness=-1
The target metric in~\eqref{MainMaps} does not have constant curvature nevertheless is globally def\/ined and
positively curved everywhere; in fact $K=2^{-1}(1+e^{u_1})^{-1}$.
However, the metric is nonetheless very special because the wave map system turns out to be Darboux integrable, as
demonstrated in~\cite{Ream}.

The Lagrangian density for this metric is
\begin{gather*}
\mathcal{L}=\frac{{u_1}_x{u_1}_y-{u_2}_x{u_2}_y}{1+e^{-u_1}},
\end{gather*}
whose Euler--Lagrange equation is
\begin{gather*}
{u_1}_{xy}+\frac{{u_1}_x{u_1}_y+{u_2}_x{u_2}_y}{2(1+e^{u_1})}=0,
\qquad
{u_2}_{xy}+\frac{{u_1}_x{u_2}_y+{u_2}_x{u_1}_y}{2(1+e^{u_1})}=0.
\end{gather*}
The change of variables $(u_1,u_2)\mapsto((u_1+u_2)/2,(u_1-u_2)/2)=(u,v)$ transforms this to
\begin{gather}
\label{Cau1}
u_{xy}+\frac{u_xu_y}{2\left(1+e^{u/2+v/2}\right)}=0,
\qquad
v_{xy}+\frac{v_xv_y}{2\left(1+e^{u/2+v/2}\right)}=0.
\end{gather}

We now prove
\begin{theorem}
\label{IVPmain}
Consider the initial value problem
\begin{gather}
\label{CauchyTheorem}
u_{xy}+\frac{u_xu_y}{2\left(1+e^{u/2+v/2}\right)}=0,
\qquad
v_{xy}+\frac{v_xv_y}{2\left(1+e^{u/2+v/2}\right)}=0,
\\
u_{|_\gamma}=\phi_1,
\qquad
v_{|_\gamma}=\phi_2,
\qquad
\frac{\partial u}{\partial\boldsymbol{n}}_{|_\gamma}=\psi_1,
\qquad
\frac{\partial v}{\partial\boldsymbol{n}}_{|_\gamma}=\psi_2,
\end{gather}
where $\gamma$ is a~curve with tangents nowhere parallel to the $x$- or $y$-axes, $\boldsymbol{n}$ is a~unit normal
vector field along $\gamma$ and $\phi_i$, $\psi_i$ are smooth functions along $\gamma$.
\begin{enumerate}\itemsep=0pt
\item[$1.$] Problem~\eqref{CauchyTheorem} has a~unique smooth local solution.
Moreover, the unique local solution is exp\-res\-sible as the solution of an ordinary differential equation of Lie type
associated an action of ${\rm SL}(2)$ on $\mathbb{R}^4$.
\item[$2.$] Given the unique local solution $(u,v)$ from part~$1$, the Cauchy problem for harmonic maps
\begin{gather*}
\big(\mathbb{R}^{1,1}, dxdy\big)\to\left(M,\frac{du_1^2-du_2^2}{1+e^{-u_1}}\right).
\end{gather*}
is given by
\begin{gather*}
u_1=u+v,
\qquad
u_2=u-v.
\end{gather*}
where $u_1$, $u_2$ satisfy initial conditions
\begin{gather*}
{u_1}_{|_\gamma}=\phi_1+\phi_2,
\qquad
{u_2}_{|_\gamma}=\phi_1-\phi_2,
\qquad
\frac{\partial u_1}{\partial\boldsymbol{n}}_{|_\gamma}=\psi_1+\psi_2,
\qquad
\frac{\partial u_2}{\partial\boldsymbol{n}}_{|_\gamma}=\psi_1-\psi_2.
\end{gather*}
\end{enumerate}
\end{theorem}

\begin{proof} By hyperbolicity, the problem is locally well-posed.
System~\eqref{Cau1} has {\it four} f\/irst integrals on each characteristic system, $H_1$, $H_2$.
Let us label the f\/irst integrals
\begin{gather*}
y,\;\beta_1,\;\beta_2,\;\beta_3
\quad
\text{for}
\quad
H_1,
\qquad
\text{and}
\qquad
x,\;\alpha_1,\;\alpha_2,\;\alpha_3
\quad
\text{for}
\quad
H_2.
\end{gather*}
For this system $\alpha_1$, $\beta_1$ are {\it first order} dif\/ferential functions while
$\alpha_2$, $\alpha_3$, $\beta_2$, $\beta_3$ are of second order.
Finally, while the 8 f\/irst integrals are functionally independent, we have
\begin{gather*}
\frac{d\alpha_1}{dx}=\alpha_2
\qquad
\text{and}
\qquad
\frac{d\beta_1}{dy}=\beta_2.
\end{gather*}

Let $k_1(y)$, $k_2(y)$ be arbitrary functions and consider the overdetermined PDE system def\/ined by~\eqref{Cau1}
together with the additional PDE
\begin{gather}
\label{Cau2}
\beta_1=k_1(y),\qquad \beta_2=\dot{k}_1(y),\qquad \beta_3=k_3(y),
\end{gather}
where the dot denotes $y$-dif\/ferentiation.
It can be shown that this overdetermined system $\mathcal{E}'$ is involutive and moreover admits a~1-dimensional Cauchy
distribution, as we will see.
Now suppose we f\/ix a~smooth curve $\gamma$ embedded in a~portion of the $xy$-plane, $\mathcal{N}$ and suppose Cauchy
data is prescribed along $\gamma$ as in the Theorem statement.
Then by an argument similar to~\cite[Proposition~3.3]{Vassiliou2000}, $\gamma$ can be lifted to a~unique curve
$\widehat{\gamma}:I\subseteq\mathbb{R}\to J^1(\mathcal{N},\mathbb{R}^2)$ which agrees with all the Cauchy data.
Let $\iota:\mathcal{H}_{k_1,k_2}\to J^2(\mathcal{N},\mathbb{R}^2)$ denote the submanifold in
$J^2(\mathcal{N},\mathbb{R}^2)$ def\/ined by PDE system $\mathcal{E}'$ and $\Theta$ the contact system on $J^2$.
Let $\theta_{k_1,k_2}=\iota^*\Theta$ be the Pfaf\/f\/ian system whose integral submanifolds are the solutions of~$\mathcal{E}'$.
Our aim is to choose the functions $k_1$, $k_2$, if possible, in order that we can extend $\widehat{\gamma}$ to
a~1-dimensional integral submanifold $\widetilde{\gamma}$ of $\theta_{k_1,k_2}$.

Now
\begin{gather*}
\theta_{k_1,k_2}=\{\omega^1=du-u_xdx-u_ydy,
\,
\omega^2=dv-v_xdx-v_ydy,
\,
\omega^3=du_x-u_{xx}dx-f^1dy,
\\
\hphantom{\theta_{k_1,k_2}=\{}{} \omega^4=du_y-f^1dx-u_{yy}dy,
\,
\omega^5=dv_x-v_{xx}dx-f^2dy,
\,
\omega^6=dv_y-f^2dx-v_{yy}dy\}.
\end{gather*}
Pulling back by $\widehat{\gamma}$ we observe that $\omega^1$, $\omega^2$ pullback to zero by construction.
Forms $\omega^4$ and $\omega^6$ def\/ine the functions $u_{yy}$ and $v_{yy}$ along $\gamma$ while $\omega^3$, $\omega^5$
def\/ine functions $u_{xx}$ and $v_{xx}$ along $\gamma$.
All these functions are expressed in terms of the Cauchy data, $\phi_i$, $\psi_i$.
Substituting these back into~\eqref{Cau2} determines the functions $k_1$, $k_2$ in terms of $\phi_i$, $\psi_i$.
For completeness we give the f\/irst integrals of~$H_2$
\begin{gather*}
x,
\qquad
\alpha_1=\frac{u_xv_x}{1+\exp\big(-\frac{u+v}{2}\big)},
\qquad
\alpha_2=\frac{d\alpha_1}{dx},
\\
\alpha_3=\frac{\big(
2v_xu_{xx}-2u_xv_{xx}-u_xv_x^2+v_xu_x^2+2\big(v_xu_{xx}-u_xv_{xx}-u_xv_x^2+u_x^2v_x\big)\exp\big(\frac{u+v}{2}\big)
\big)}{u_xv_x\big(1+\exp\big(\frac{u+v}{2}\big)\big)}.
\end{gather*}
Those of $H_1$ are similar but with $y$ replacing $x$.

We will now use these f\/irst integrals to demonstrate that $\theta_{k_1,k_2}$ has a~one-dimensional Cauchy distribution
and that in particular the Cauchy vector can be chosen to be a~curve in a~certain Lie algebra~-- the Vessiot
algebra~\cite{AFV} of system~\eqref{Cau1}.
It will be seen that the Cauchy vector is generically transverse to the Cauchy data and extends the one-dimensional
integral submanifold of $\theta_{k_1,k_2}$ to the solution of the Cauchy problem.
Because the Cauchy vector is a~curve in a~Lie algebra, this extension from a~one-dimensional to a~two-dimensional
integral of $\theta_{k_1,k_2}$ is an ordinary dif\/ferential equation~$\mathfrak{L}$ of Lie type.
Its coef\/f\/icients and initial conditions are f\/ixed by all the data present in the problem, including the Cauchy
data.
{\it Any} solution of $\mathfrak{L}$ (independently of its initial conditions) permits a~Lie reduction of
$\mathfrak{L}$ and will permit us to solve the IVP for $\mathfrak{L}$.

Indeed, setting $z_1=e^{u/2}$, $z_2=e^{v/2}$, $z_3=u_y$, $z_4=v_x$, $a_i=\alpha_i$, $b_i=\beta_i$, $i=1,2,3$ we calculate that
the contact system on $J^2(\mathbb{R}^2,\mathbb{R}^2)$ pulled back to PDE~\eqref{Cau1} is
$\Psi=\{\kappa^1,\ldots,\kappa^6\}$, where
\begin{gather*}
\kappa^1=da_1-a_2dx,
\qquad
\kappa^2=db_1-b_2dy,
\qquad
\kappa^3=dz_1-\frac{a_1(1+z_1z_2)}{2z_2z_4}dx-\frac{z_1z_3}{2}dy,
\\
\kappa^4=
dz_2-\frac{z_2z_4}{2}dx-\frac{b_1(1+z_1z_2)}{2z_1z_3}dy,
\qquad
\kappa^5=dz_3+\frac{a_1z_3}{2z_1z_2z_4}dx+\frac{1}{2}(z_3^2-b_3z_3-b_1)dy,
\\
\kappa^6=dz_4+\frac{1}{2}(z_4^2-a_3z_4-a_1)dx+\frac{b_1z_4}{2z_1z_2z_3}dy.
\end{gather*}
Pulling $\Psi$ back to submanifold~\eqref{Cau2} yields a~Pfaf\/f\/ian system with 1-dimensional Cauchy distribution
spanned by
\begin{gather*}
\xi_{k_1,k_2}=\P y-\frac{k_2(y)}{4}\left(R_1+4R_4\right)+\frac{k_1(y)}{2}R_2+\frac{1}{2}R_3,
\end{gather*}
where the $R_i$ form a~basis for the Vessiot algebra\footnote{A brief geometric construction and interpretation of the
Vessiot algebra is given in Section~\ref{section4} and~Appendix~\ref{appendix}.
See~\cite{AFV} for a~complete exposition.
However, detailed knowledge of Vessiot algebras is not a~prerequisite for this paper.}
\begin{gather*}
\begin{split}
& R_1=z_1\P{z_1}-z_2\P{z_2}-2z_3\P{z_3},
\qquad
R_2=\frac{1+z_1z_2}{z_1z_3}\P{z_2}+\P{z_3}-\frac{z_4}{z_1z_2z_3}\P{z_4},
\\
& R_3=z_1z_3\P{z_1}-z_3^2\P{z_3},
\qquad
R_4=-\frac{1}{4}\left(z_1\P{z_1}-z_2\P{z_2}\right)
\end{split}
\end{gather*}
with nonzero structure
\begin{gather*}
[R_1,R_2]=2R_2,
\qquad
[R_1,R_3]=-2R_3,
\qquad
[R_2,R_3]=R_1.
\end{gather*}
Since $\xi_{k_1,k_2}$ is a~curve in the Vessiot algebra it determines an ODE of Lie type.
Furthermore $\xi_{k_1,k_2}$ is generically transverse to the Cauchy data.

Note that $\rho_1=R_1+4R_4$, $\rho_2=R_2$, $\rho_3=R_3$ generates a~local action of ${\rm SL}(2)$ on $\mathbb{R}^4$:
\begin{gather*}
[\rho_1,\rho_2]=2\rho_2,
\qquad
[\rho_1,\rho_3]=-2\rho_3,
\qquad
[\rho_2,\rho_3]=\rho_1
\end{gather*}
and the Cauchy vector is
\begin{gather*}
\xi_{k_1,k_2}=\P y-k_2(y)\rho_1+k_1(y)\rho_2+\frac{1}{2}\rho_3.
\end{gather*}
In summary, vector f\/ield $\xi_{k_1,k_2}$ f\/lows the 1-dimensional initial data solution curve $\widetilde{\gamma}$
of $\theta_{k_1,k_2}$ to a~2-dimensional solution.
This completes the proof of Theorem~\ref{IVPmain}.
\end{proof}

\begin{example}\label{example1}
As an illustrative example we consider the initial value problem
\begin{gather*}
u_{xy}+\frac{u_xu_y}{2\left(1+e^{u/2+v/2}\right)}=0,
\qquad
v_{xy}+\frac{v_xv_y}{2\left(1+e^{u/2+v/2}\right)}=0,
\\
u_{|_\gamma}=v_{|_\gamma}=0,
\qquad
\frac{\partial u}{\partial\boldsymbol{n}}_{|_\gamma}=\frac{\partial v}{\partial\boldsymbol{n}}_{|_\gamma}=\sqrt{2},
\end{gather*}
where $\gamma=(x,x)$.
Since~$x$ and~$y$ are light-cone coordinates $x=(\xi+\tau)/2$, $y=(\xi-\tau)/2$, the curve $\gamma$ corresponds to time
$\tau=0$.
Thus we have constant initial values at time $\tau=0$.
We wish to determine the system of Lie type whose solutions corresponds to the solution of this Cauchy problem.
We have $\boldsymbol{n}={2}^{-1/2}(\P x-\P y)$ and we get
\begin{gather*}
{u_x}_{|_\gamma}={v_x}_{|_\gamma}=1,
\qquad
{u_y}_{|_\gamma}={v_y}_{|_\gamma}=-1.
\end{gather*}
\end{example}

So our initial curve in $J^1$ is $\widehat{\gamma}(x)=(x,x,0,0,1,-1,1,-1)=(x,y,u,v,u_x,u_y,v_x,v_y)$.
This translates to an initial curve in the adapted coordinates
\begin{gather*}
(x,y,z_1,z_2,z_3,z_4)=(x,x,1,1,-1,1).
\end{gather*}
We extend this to a~unique 1-dimensional integral of $\Psi$ and get
\begin{gather*}
\widetilde{\gamma}(x)=(x,y,z_1,z_2,z_3,z_4,a_1,a_2,a_3,b_1,b_2,b_3)=
\left(x,x,1,1,-1,1,\frac{1}{2},0,0,\frac{1}{2},0,0\right).
\end{gather*}
Thus we get
\begin{gather*}
k_1(y)=1/2,\qquad k_2(y)=0.
\end{gather*}
The Cauchy vector is therefore
\begin{gather*}
\xi_{\frac{1}{2},0}=\P y+\frac{1}{4}R_2+\frac{1}{2}R_3.
\end{gather*}
We f\/low this vector f\/ield obtaining solutions $z_i(x,y)$ subject to the intial conditions
\begin{gather*}
z_1(x,x)=1,\qquad z_2(x,x)=1,\qquad z_3(x,x)=-1,\qquad z_4(x,x)=1.
\end{gather*}

The ODE are
\begin{gather*}
\frac{\partial z_1}{\partial y}=z_1z_2,
\qquad
\frac{\partial z_2}{\partial y}=\frac{1+z_1z_2}{4z_1z_3},
\qquad
\frac{\partial z_3}{\partial y}=\frac{1}{4}-\frac{z_3^2}{2},
\qquad
\frac{\partial z_4}{\partial y}=\frac{z_4}{4z_1z_2z_3}
\end{gather*}
to be solved for $z_i(x,y)$ subject to the given initial conditions.
In fact, we need not solve the whole system but only the equation for $z_3$ and then substitute this into the equation
for~$z_1$.
This gives the function $u$ up to a~quadrature after which the function $v$ can be obtained algebraically from the PDE
itself
\begin{gather*}
e^{v/2}=-e^{-u/2}\left(1+\frac{u_xu_y}{2u_{xy}}\right).
\end{gather*}
We get the unique solution of the Cauchy problem to be
\begin{gather*}
u(x,y)=v(x,y)=\ln\left(\cosh\frac{\sqrt{2}}{4}(x-y)+\sqrt{2}\sinh\frac{\sqrt{2}}{4}(x-y)\right)^2.
\end{gather*}

\begin{example}\label{example2} A slightly more interesting example is obtained from the initial conditions
\begin{gather*}
u_{|_\gamma}=2\ln\lambda,\qquad v_{|_\gamma}=2\ln\frac{1}{\lambda},\qquad \frac{\partial u}{\partial\boldsymbol{n}}_{|_\gamma}=
\frac{\partial v}{\partial\boldsymbol{n}}_{|_\gamma}=1,
\end{gather*}
for any $\lambda>0$.
The Cauchy vector is $\xi_{\frac{1}{4},0}$ and its f\/low subject to the initial conditions along $y=x$ being
\begin{gather*}
z_1=\lambda,\qquad z_2=\frac{1}{\lambda},\qquad z_3=-\frac{1}{\sqrt{2}},\qquad z_4=\frac{1}{\sqrt{2}}
\end{gather*}
gives rise to the unique solution
\begin{gather*}
u=2\ln\left(\lambda\exp\left(\frac{x+y}{4}\right)\frac{\exp\left(-\frac{y}{2}\right)(3+2\sqrt{2})-\exp\left(-\frac{x}{2}\right)}{2(1+\sqrt{2})}\right),
\\
v=2\ln\left(\frac{1}{2\lambda}\exp\left(\frac{x+y}{4}\right)\frac{\exp\left(-y\right)\left(2\sqrt{2}+3\right)-\exp\left(-x\right)\left(2\sqrt{2}-3\right)-2\exp\left(-\frac{x+y}{2}\right)}{\exp\left(-\frac{y}{2}\right)\left(1+\sqrt{2}\right)+\exp\left(-\frac{x}{2}\right)\left(1-\sqrt{2}\right)}\right).
\end{gather*}
Thus, even constant initial data has the potential of producing some interesting explicit solutions.
Indeed, one can ask if this solution is global in time.
\end{example}

The fact that we only had to solve for~$z_3$ in Example~\ref{example1} (and Example~\ref{example2}) holds not only for
these illustrative examples since the system of Lie type in the general case is
\begin{gather*}
\frac{\partial z_1}{\partial y}=\frac{1}{2}z_1z_3,
\qquad
\frac{\partial z_2}{\partial y}=k_1(y)\frac{(1+z_1z_2)}{2z_1z_3},
\\
\frac{\partial z_3}{\partial y}=\frac{1}{2}\big(k_1(y)+k_2(y)z_3-z_3^2\big),
\qquad
\frac{\partial z_4}{\partial y}=-\frac{k_1(y)z_4}{2z_1z_2z_3}.
\end{gather*}
This proves

\begin{theorem}
\label{RiccatiReductionThm}
In the Cauchy problem for wave maps
\begin{gather*}
\big(\mathbb{R}^{1,1}, dx dy\big)\to\left(M,\frac{du_1^2-du_2^2}{1+e^{-u_1}}\right)
\end{gather*}
of Theorem~{\rm \ref{IVPmain}}, let functions $k_1(y)$ and $k_2(y)$ be defined as described above.
Denote by~$\zeta$ the initial value of $z_3=u_y$ along $\gamma$.
Let~$\Gamma$ be the unique function satisfying Riccati initial value problem
\begin{gather}
\label{RiccatiIVP}
\frac{\partial\Gamma}{\partial y}=\frac{1}{2}\big(k_1(y)+k_2(y)\Gamma-\Gamma^2\big),
\qquad
\Gamma_{|_\gamma}=\zeta.
\end{gather}
Define function $u$ to be the unique solution of
\begin{gather*}
\frac{\partial u}{\partial y}=\Gamma,
\qquad
{u}_{|_\gamma}=\phi_1
\end{gather*}
and let $v$ be defined algebraically from the partial differential equation
\begin{gather*}
u_{xy}+\frac{u_xu_y}{2\left(1+e^{u/2+v/2}\right)}=0
\end{gather*}
upon substituting solution $u$ and solving for $v$.
The functions $u$, $v$ constitute the unique solution of the Cauchy problem.
\end{theorem}

This implies that the solution of any given Cauchy problem for wave maps into the metric
\begin{gather*}
h=\frac{du_1^2-du_2^2}{1+e^{-u_1}}
\end{gather*}
relies on the solution of a~Riccati initial value problem together with one quadrature.
The interesting point here is that the Riccati equation is the {\it simplest} non-trivial equation of Lie type.
It is a~Lie system for the local ${\rm SL}(2)$-action on the real line that globalises on $\mathbb{RP}^1$.
In general, solutions of Riccati equations develop singularities in f\/inite time, even those with constant
coef\/f\/icients.
However, the theorem above provides a~correspondence between Cauchy data for the wave map and the Riccati initial value
problem~\eqref{RiccatiIVP}.
An interesting problem is to study this correspondence and link the nature of the Cauchy data with the properties of
the solution of~\eqref{RiccatiIVP} and in turn, link this correspondence with the geometry of the target metric.

For the standard initial value problem where Cauchy data is posed along $0=2\tau=x-y$, the relationship between the
Cauchy data and the coef\/f\/icients of the Riccati equation is complicated.
However, due to its signif\/icance and for latter use, we give it explicitly:
\begin{gather*}
k_1=
-\frac{1}{4}\frac{\psi_2\left(\sqrt{2}{\phi_1}_x-2\psi_1\right)\exp\left(\frac{\phi_2}{2}\right)}
{\exp\left(-\frac{\phi_1}{2}\right)+\exp\left(\frac{\phi_2}{2}\right)},
\qquad
a_1=\frac{1}{4}\frac{(\sqrt{2}{\phi_1}_x+2\psi_1)(\sqrt{2}{\phi_2}_x+\psi_2)
\exp\left(\frac{\phi_1}{2}\right)}{\exp\left(\frac{\phi_1}{2}\right)+\exp\left(-\frac{\phi_2}{2}\right)},
\\
k_2=
\delta^{-1}\Bigg(2\left({\phi_1}^2_x+2\psi_2^2-2\sqrt{2}\psi_1{\phi_1}_x-4k_1\right)
\left(\psi_2+\sqrt{2}{\phi_2}_x\right)
\\
\phantom{k_2=}{}
+4\sqrt{2}a_1\left({\phi_1}_x-\psi_1\sqrt{2}\right)
\exp\left(-\frac{\phi_1}{2}-\frac{\phi_2}{2}\right)+8\left(-{\psi_1}_x\sqrt{2}+{\phi_1}_{xx}\right)
\left(\psi_2+\sqrt{2}{\phi_2}_x\right)\Bigg),
\end{gather*}
where
\begin{gather*}
\delta=4\left(\sqrt{2}{\phi_2}_x+\psi_2\right)\left({\phi_1}_x-\psi_1\sqrt{2}\right).
\end{gather*}
Recall that functions $\phi_1$, $\phi_2$ are respectively the values of $u$ and $v$ along $\gamma$, the $\psi_i$ are the
values of the normal derivatives along $\gamma$ as stated in the theorem.
For instance for arbitrary constant initial conditions along the time axis
\begin{gather*}
\phi_1(x)=u(x,x)=2\ln\lambda_1,
\qquad
\phi_2(x)=v(x,x)=2\ln\lambda_2,
\\
\psi_1=\alpha\sqrt{2},
\qquad
\psi_2=\beta\sqrt{2};
\qquad
\lambda_1,\; \lambda_2>0,
\end{gather*}
we have
\begin{gather}
\label{constantIC}
k_1(y)=\frac{\alpha\beta\lambda_1\lambda_2}{1+\lambda_1\lambda_2},
\qquad
k_2(y)=-\frac{(\alpha-\beta)\lambda_1\lambda_2}{1+\lambda_1\lambda_2}.
\end{gather}
The solution of the Cauchy problem in this case is easily calculated as we did in Examples~\ref{example1} and~\ref{example2} because it amounts
to solving Riccati equation~\eqref{RiccatiIVP} with constant coef\/f\/icients $k_1(y)$, $k_2(y)$ given
by~\eqref{constantIC}; indeed Example 1 is the choice $\lambda_1=\lambda_2=1$, $\alpha=\beta=\sqrt{2}$.
However the general formula is complicated and of itself not very informative so we refrain from recording it here.

\section{Cauchy problem for wave maps and Lie reduction} \label{section3}

We have shown that to solve the Cauchy problem for the wave
map equation it is enough to ``evolve" the initial data curve $\widetilde{\gamma}$ by solving a~Riccati equation which
is notably the simplest nontrivial system of Lie type.
One signif\/icant feature of Lie systems is that they admit ``reduction by particular solutions'', otherwise known as
{\it Lie reduction}.
There are only a~few sources scattered in the literature on this topic, among
them~\cite{Bryant,Carinena,DoubrovKomrakov}.
In this section we give an illustration of how Lie reduction may be useful in resolving instances of the Cauchy problem
for our wave map system.
Appendix~\ref{appendixA2} of this paper summarises the known results on systems of Lie type, oriented toward the
applications at hand.
In this and subsequent sections we will adopt the notation and theory set out in Appendix~\ref{appendixA2}, to which we refer the
reader.

Consider the Cauchy problem for wave maps with (non-constant) Cauchy data
\begin{gather*}
\phi_1=\phi_2=0,\qquad \psi_1=-2\sqrt{2},\qquad \psi_2=2\sqrt{2}x
\end{gather*}
along the curve $y=x$.
We get
\begin{gather*}
k_1(y)=-2y,
\qquad
k_2(y)=-y-1,
\end{gather*}
so that the corresponding Riccati initial value problem of Theorem~\ref{RiccatiReductionThm} is
\begin{gather}
\label{RiccatiIVPexmple}
\frac{\partial\Gamma}{\partial y}=\frac{1}{2}\big({-}2y-(y+1)\Gamma-\Gamma^2\big),
\qquad
 \Gamma(x,x)=2.
\end{gather}
We observe that $\Gamma=1-y$ is a~solution vanishing at $y=1$.
Implementing the procedure described in Appendix~\ref{appendixA2} obtains one factor in the fundamental solution
\begin{gather*}
g_0(y)=\left(\begin{matrix}1&1-y\\0&1\end{matrix}\right).
\end{gather*}
A curve in the isotropy group of 0 is denoted by $H$ and has the form
\begin{gather*}
g_1(y)=\left(\begin{matrix}\gamma_1(y)&0\\\gamma_2(y)&\gamma_1(y)^{-1}\end{matrix}\right).
\end{gather*}
The curve of Lie algebra elements associated to the Riccati equation
\begin{gather*}
\frac{dz}{dy}=\alpha_0(y)+2\alpha_1(y)z-\alpha_2(y)z^2
\end{gather*}
is
\begin{gather*}
A(y)=\left(\begin{matrix}\alpha_1&\alpha_0\\\alpha_2&-\alpha_1\end{matrix}\right)\subset\mathfrak{sl}(2).
\end{gather*}
For the Riccati initial value problem~\eqref{RiccatiIVPexmple} we have
\begin{gather*}
A(y)=\left(\begin{matrix}-\frac{y+1}{4}&y\\\frac{1}{2}&\frac{y+1}{4}\end{matrix}\right)
\end{gather*}
and the reduced fundamental equation is
\begin{gather*}
\frac{dg_1}{dy}=B(y)g_1,
\end{gather*}
where (see Appendix~\ref{appendixA2}, especially Theorem~\ref{LieRed} \& Appendix~\ref{LieRedRiccati})
\begin{gather*}
B(y)=g_0(y)^{-1}A(y)g_0(y)-g_0(y)^{-1}\frac{dg_0}{dy}=
\left(\begin{matrix}\frac{1}{4}y-\frac{3}{4}&0\\\frac{1}{2}&-\frac{1}{4}y+\frac{3}{4}\end{matrix}\right),
\end{gather*}
valued in the isotropy subalgebra at $0$, as expected.
The ODE initial value problem for the fundamental solution factor $g_1(y)$ is
\begin{gather*}
\frac{d\gamma_1}{dy}-\beta_1\gamma_1=0,
\qquad
\frac{d\gamma_2}{dy}-\frac{1}{2}\gamma_1+\beta_1\gamma_2=0,
\qquad
\gamma_1(1)=1,\qquad \gamma_2(1)=0,
\end{gather*}
where $\beta_1=(y-3)/4$.
This problem can be explicitly solved in terms of elementary functions giving
\begin{gather*}
\gamma_1(y)=\exp\left(\frac{(y-1)(y-5)}{8}\right),
\\
\gamma_2(y)=\sqrt{\pi}\exp\left(\frac{3}{4}y-y^2-\frac{13}{8}\right)
\left(\text{erf\/i}(1)+\text{erf\/i}\left(\frac{1}{2}\left(y-3\right)\right)\right),
\end{gather*}
where $i=\sqrt{-1}$ and $\text{erf\/i}$ denotes a~concomitant of the error function:
\begin{gather*}
\text{erf\/i}(x)=\frac{2}{\sqrt{\pi}}\int_0^x\exp\big(t^2\big) \,dt.
\end{gather*}
This data enables us to construct the fundamental solution $g(y)=g_0(y)g_1(y)$ for~\eqref{RiccatiIVPexmple} and leads
to the solution of the ODE in~\eqref{RiccatiIVPexmple}
\begin{gather*}
\Gamma(y;q)=\lambda_{g(y)}(q),
\end{gather*}
where $\lambda_h$ denotes the linear fractional transformation~\eqref{Mobius} by $h\in {\rm SL}(2)$.

It is now a~simple matter to determine the value of $q$ that satisf\/ies the initial condition $\Gamma(x,x)=2$ and
giving the unique solution $\Gamma(x,y)$ of the Riccati initial value problem~\eqref{RiccatiIVPexmple}.
We f\/ind
\begin{gather*}
\Gamma(x,y)=
1-y-\frac{2\exp\left(\frac{3}{2}y-\frac{1}{4}y^2\right)}{\Delta(x)-\sqrt{\pi}\exp\left(-\frac{9}{4}\right)
\text{erf\/i}\left(\frac{1}{2}(y-3)\right)},
\end{gather*}
where
\begin{gather*}
\Delta(x)=
-\frac{2\exp\left(-\frac{3}{2}x+\frac{1}{4}x^2\right)}{(1+x)}+\sqrt{\pi}\exp\left(-\frac{9}{4}\right)\text{erf\/i}\left(\frac{1}{2}(x-3)\right).
\end{gather*}
Finally, we obtain an integral representation of the solution $u$ satisfying
\begin{gather*}
\frac{\partial u}{\partial y}=\Gamma(x,y),
\qquad
u(x,x)=0,
\end{gather*}
namely
\begin{gather*}
u(x,y)=\int_{x}^y\Gamma(x,s)\,ds
\end{gather*}
or in terms of spacetime coordinates $(\xi,\tau)$
\begin{gather*}
\bar{u}(\xi,\tau)=\int_{\frac{1}{2}(\xi+\tau)}^{\frac{1}{2}(\xi-\tau)}\Gamma\left(\frac{1}{2}(\xi+\tau),s\right) ds.
\end{gather*}
The signif\/icance of Lie reduction in our ability to solve the Riccati equation should here be emphasised.
Without this, it would have been impossible to construct the fundamental solution and there would be no hope of
constructing function $\Gamma$ and constructing the integral representation of the solution of a~Cauchy problem with
non-constant initial data.

In the example above we relied on knowledge of a~simple solution, namely $\Gamma(x,y)=1-y$ to perform the reduction.
But even with polynomial or rational coef\/f\/icients a~Riccati equation will not generally have any rational solutions.
In this case however we can appeal to the well known fact that every Riccati equation can be linearised.
\begin{lemma}
\label{RiccatiLinearization}
The general Riccati equation
\begin{gather}\label{gRe}
\frac{dz}{dt}=\alpha_0(t)+2\alpha_1(t)z(t)-\alpha_2(t)z(t)^2
\end{gather}
can be transformed to the form
\begin{gather*}
\frac{d}{d\tau}y(\tau)=\beta(\tau)+y(\tau)^2,
\end{gather*}
where $z=p(t)y(t)$
\begin{gather*}
\ln{p}=\int^t2\alpha_1(s)\,ds,
\qquad
\tau=\int^t\alpha_2(s)p(s)\,ds.
\end{gather*}
Provided these quadratures can be carried out then the explicit solvability of~\eqref{gRe} depends on the
properties of its linearisation
\begin{gather*}
\frac{d^2\psi}{d\tau^2}+\beta(\tau)\psi=0,
\qquad
\text{where}
\qquad
 y(\tau)=-\frac{1}{\psi(\tau)}\frac{d}{d\tau}\psi(\tau).
\end{gather*}
Any solution of the $2$nd order linear ODE can be used in the Lie reduction of the Riccati equation.
\end{lemma}

As a~consequence of Lemma~\ref{RiccatiLinearization} and Theorem~\ref{quadReduction} of Appendix~\ref{appendixA2}, we have
\begin{corollary} The solution of the Cauchy problem for wave maps of Theorem~{\rm \ref{IVPmain}}
\begin{gather*}
\big(\mathbb{R}^{1,1}, dx dy\big)\to\left(M,\frac{du_1^2-du_2^2}{1+e^{-u_1}}\right)
\end{gather*}
is reducible to quadrature provided a~particular solution of the Riccati equation~\eqref{RiccatiIVP} is known.
A~particular solution of~\eqref{RiccatiIVP} can be constructed by quadrature and the solution of a~linear second order
ODE.
\end{corollary}

\begin{remark} As a~consequence of Lemma~\ref{RiccatiLinearization} and Theorem~\ref{quadReduction} the dif\/ferential
equations solver in {\tt MAPLE}~-- {\tt dsolve } is very often able to construct an explicit representation in terms of
known special functions to a~Riccati initial value problem when the coef\/f\/icients are polynomial functions of the
independent variable.
\end{remark}

\section{Hyperbolic Weierstrass representation}\label{section4}

In this section we use the Darboux integrability of the wave map equation~\eqref{Cau1} to compute its general solution
and hence construct a~hyperbolic Weierstrass-type representation for wave maps into the corresponding metric.
According to~\cite{AFV}, we pull back $\Psi$ to a~suitable integral manifold $M_1,M_2$ of $H_1^{(\infty)}$ and
$H_2^{(\infty)}$ respectively.
It is convenient to def\/ine $M_1$ by $y=b_1=b_2=b_3=0$ and $M_2$ by $x=a_1=a_2=a_3=0$.
This gives Pfaf\/f\/ian systems
\begin{gather*}
\Psi_1=\bigg\{dz_4-\frac{1}{2}\left(a_3z_4+a_1-z_4^2\right)dx,\,
dz_3+\frac{1}{2z_1z_2z_4}a_1z_3dx,\,
dz_2-\frac{1}{2}z_2z_4dx,
\\
\phantom{\Psi_1=\bigg\{}
dz_1-\frac{1}{2z_2z_4}a_1(1+z_1z_2)dx,\,
da_1-a_2dx\bigg\}
\end{gather*}
and
\begin{gather*}
\Psi_2=\bigg\{dz_1-\frac{1}{2}z_3z_1dy,\,
dz_2-\frac{1}{2z_1z_3}b_1(1+z_1z_2)dy,\,
dz_4+\frac{1}{2z_1z_2z_3}b_1z_4dy,
\\
\phantom{\Psi_1=\bigg\{}
dz_3-\frac{1}{2}\big(b_3z_3-z_3^2+b_1\big)dy,\,
db_1-b_2dy\bigg\},
\end{gather*}
each of rank 5 on 8-manifolds: $(M_i$, $\Psi_i)$.
Locally $M_1$ has coordinates $x$, $z_1$, $z_1$, $z_2$, $z_3$, $z_4$, $a_1$, $a_2$, $a_3$ while $M_2$ has local coordinates
$y$, $z_1$, $z_1$, $z_2$, $z_3$, $z_4$, $b_1$, $b_2$, $b_3$.
Using these local formulas, we def\/ine a~local product structure
\begin{gather*}
\left(\widehat{M}_1\times\widehat{M}_2,\widehat{\Psi}_1\oplus\widehat{\Psi}_2\right),
\end{gather*}
where
\begin{gather*}
\widehat{\Psi}_1=\bigg\{dq_4-\frac{1}{2}\left(a_3q_4+a_1-q_4^2\right)dx,
\,
dq_3+\frac{1}{2q_1q_2q_4}a_1q_3dx,
\,
dq_2-\frac{1}{2}q_2q_4dx,
\\
\phantom{\widehat{\Psi}_1=\bigg\{}
dq_1-\frac{1}{2q_2q_4}a_1(1+q_1q_2)dx,
\,
da_1-a_2dx\bigg\}
\end{gather*}
and
\begin{gather*}
\widehat{\Psi}_2=\bigg\{dp_1-\frac{1}{2}p_1p_3dy,
\,
dp_2-\frac{1}{2p_1p_3}b_1(1+p_1p_2)dy,
\,
dp_4+\frac{1}{2p_1p_2p_3}b_1p_4dy,
\\
\phantom{\widehat{\Psi}_2=\bigg\{}
dp_3-\frac{1}{2}\big(b_3p_3-p_3^2+b_1\big)dy,
\,
db_1-b_2dy\bigg\}.
\end{gather*}

As described in~\cite{AFV}, the relationship between $\Psi$, $\widehat{\Psi}_1$ and $\widehat{\Psi}_2$ is that every
integral manifold of $\Psi$ is a~superposition of an integral manifold of $\widehat{\Psi}_1$ and $\widehat{\Psi}_2$.
The superposition formula is the map\footnote{See Appendix~\ref{appendix} for further details on the superposition
formula; in particular, how it is def\/ined and constructed.}
\begin{gather*}
\boldsymbol{\sigma}: \ \widehat{M}_1\times\widehat{M}_2\to M
\end{gather*}
def\/ined by
\begin{gather*}
\boldsymbol{\sigma}(x,\,\boldsymbol{p},\,\boldsymbol{a};\,y,\,\boldsymbol{q},\,\boldsymbol{b})=
\Bigg(\frac{1-p_1p_2p_4+2p_1p_2p_4q_3-p_1p_2q_3-q_3+p_1p_2)q_1}{p_4p_2},
\\
\qquad\quad
\frac{p_2(2p_4q_3q_2q_1-q_1q_2q_3-p_4q_2q_1+q_1q_2-p_4+1)}{q_3q_1},
\\
\qquad\quad
\frac{(2p_4q_3-q_3-p_4+1)p_3p_1p_2}{1-p_1p_2p_4+2p_1p_2p_4q_3-p_1p_2q_3-q_3+p_1p_2},
\\
\qquad\quad
\frac{(2p_4q_3-q_3-p_4+1)q_1q_2q_4}{2p_4q_3q_2q_1-q_1q_2q_3-p_4q_2q_1+q_1q_2-p_4+1},
\,x,\,y,\,a_1,\,a_2,\,a_3,\,b_1,\,b_2,\,b_3\Bigg)
\\
\qquad
=(z_1,\,z_2,\,z_3,\,z_4,\,x,\,y,\,a_1,\,a_2,\,a_3,\,b_1,\,b_2,\,b_3).
\end{gather*}
The usefulness of this factorisation of the integration problem for~$\Psi$ is not only that the integration of~$\Psi_i$
relies on ODE while that of $\Psi$ relies on PDE but that the~$\Psi_i$ are locally equivalent to prolongations of the
contact system on $J^1(\mathbb{R},\mathbb{R}^2)$.
To see this we turn to the characterisation of partial prolongations of such contact systems provided
by~\cite{Vassiliou2006a,Vassiliou2006b} which also provide simple procedures for f\/inding the equivalence.
To implement this we compute the annihilators
\begin{gather*}
\text{ann}\,\widehat{\Psi}_1:=\widehat{H}_1
\\
=\Bigg\{\P x\!+a_2\P{a_1}\!+\frac{a_1(1+q_1q_2)}{2q_2q_4}\P{q_1}\!+\frac{q_2q_4}{2}\P{q_2}\!
-\frac{a_1q_3}{2q_1q_2q_4}\P{q_3}\!+\frac{1}{2}\big(a_3q_4\!+a_1\!-q_4^2\big)\P{q_4},\, \P{a_2},\,\P{a_3}\!\Bigg\},
\\
\text{ann}\,\widehat{\Psi}_2:=\widehat{H}_2
\\
=\Bigg\{\P y\!+b_2\P{b_1}\!+\frac{p_1p_3}{2}\P{p_1}\!+\frac{b_1(1+p_1p_2)}{2p_1p_3}\P{p_2}\!
+\frac{1}{2}\big(b_3p_3-p_3^2+b_1\big)\P{p_3}\!-\frac{b_1p_4}{2p_1p_2p_3}\P{p_4},\, \P{b_2},\P{b_3}\!\Bigg\}.
\end{gather*}
We show that each of the $\widehat{H}_i$ is locally equivalent to the partial prolongation $C\langle0,1,1\rangle$ of the
contact distribution on $J^1(\mathbb{R},\mathbb{R}^2)$.
That is, the contact distribution on $J^1(\mathbb{R},\mathbb{R}^2)$ partially prolonged so that one dependent variable
has order 2 and the other order 3 with canonical local normal form
\begin{gather*}
C\langle0,1,1\rangle=\big\{\P t+z^1_1\P{z^1}+z^1_2\P{z^1_1}+z^2_1\P{z^2}+z^2_2\P{z^2_1}+z^2_3\P{z^2_2},
\,
\P{z^1_2},
\,
\P{z^2_3}\big\}.
\end{gather*}
Let $\mathcal{D}$ be a~smooth distribution on manifold $M$ and assume that $\mathcal{D}$ is totally regular in the
sense that $\mathcal{D}$, all its derived bundles and all their corresponding Cauchy bundles have constant rank on~$M$.
Denote $m_i=\dim\mathcal{D}^{(i)}$, $\chi^i=\dim \ch\mathcal{D}^{(i)}$ and
$\chi^i_{i-1}=\dim \ch\mathcal{D}^{(i)}_{i-1}$, where
\begin{gather*}
\ch\mathcal{D}^{(i)}_{i-1}=\mathcal{D}^{(i-1)}\cap\ch\mathcal{D}^{(i)}.
\end{gather*}
Below, $k$ denotes the derived length of $\mathcal{D}$.

According to Theorem~4.1 of~\cite{Vassiliou2006a}, a~totally regular distribution~$\mathcal{D}$ on smooth
manifold $M$ is locally equivalent to a~partial prolongation of the contact distribution on
$J^1(\mathbb{R},\mathbb{R}^q)$ for some~$q$ if and only if
\begin{enumerate}\itemsep=0pt
\item The integers $m_i$, $\chi^j$, $\chi^j_{j-1}$ satisfy the numerical constraints
\begin{gather*}
\chi^j=2m_j-m_{j+1}-1,
\quad
0\leq j\leq k-1,
\qquad
\chi^i_{i-1}=m_i-1,
\quad
1\leq i\leq k-1.
\end{gather*}
\item If $m_k-m_{k-1}>1$ then a~certain canonically associated bundle called the {\it resolvent} is
integrable\footnote{In the original formulation of Theorem 4.1 in~\cite{Vassiliou2006a}, the integrability of
$\ch\mathcal{D}^{(i)}_{i-1}$ is an additional hypothesis to be checked.
This is a~simple task but unnecessary since it is easy to see that this bundle is always integrable and the hypothesis
can be omitted.}.
\end{enumerate}

A pair $(M,\mathcal{D})$ that satisf\/ies these conditions is said to be a~{\it Goursat manifold} or {\it Goursat
bundle}.
Moreover, if $\mathcal{D}$ is a~Goursat bundle on $M$ then it is locally equivalent to a~partial prolongation with
$\chi^j-\chi^j_{j-1}$ dependent variables at order $j<k$ and $m_k-m_{k-1}$ dependent variables at highest order,~$k$.
This uniquely identif\/ies the partial prolongation associated to a~given Goursat manifold.
Before discussing the question of {\it constructing} equivalences, let us solve the recognition problem for the
distributions $\widehat{H}_i$.
We take $\widehat{H}_1$ which we denote temporarily by $\widehat{K}$.
We f\/ind
\begin{gather*}
\ch\widehat{K}=\{0\},
\qquad
\ch\widehat{K}^{(1)}=\{\P{a_2},\,\P{a_3}\},
\qquad
\ch\widehat{K}^{(2)}_1=\{\P{a_2},\,\P{a_3},\,\P{a_1},\,\P{q_3}\},
\\
\ch\widehat{K}^{(2)}=\{\P{a_1},\,\P{a_2},\,\P{a_3},\,\P{q_3},\,\P x\}.
\end{gather*}
Calculation shows that the dimensions of the derived bundles are
\begin{gather*}
\dim \widehat{K}=3,\qquad \dim \widehat{K}^{(1)}=5,\qquad \dim \widehat{K}^{(2)}=7,\qquad \dim \widehat{K}^{(s)}=8,
\quad
s\geq3.
\end{gather*}
Hence the derived length is $k=3$.
Below we check the f\/irst condition of a~Goursat bundle.
\begin{center}
{\bf Table.} Checking the numerical constraints satisf\/ied by $(M,\mathcal{D})$.
\vspace{1mm}

\begin{tabular}
{|c|c|c|c|c|c|}
\hline
$j$&$m_j$&$m_{j-1}-1$&$2m_j-m_{j+1}-1$&$\chi^j_{j-1}$&$\chi^j$\tsep{2pt}\bsep{1pt}\\
\hline
0&3&$-$&$6-5-1=0$&$-$&0\\ 1&5&2&$10-7-1=2$&2&2\\ 2&7&4&$14-8-1=5$&$4$&5\\
\hline
\end{tabular}
\end{center}

Hence, $\widehat{K}$ is a~Goursat bundle with $k=3$, $m_k-m_{k-1}=1$ and the only nonzero dif\/ference
$\chi^j-\chi^j_{j-1}$ being at order $j=2$: $\chi^2-\chi^2_1=5-4=1$.
Hence there is one variable of order~2 and one variable of order~3.
This solves the recognition problem and we can assert that $\widehat{K}:=\widehat{H}_1$ is locally equivalent to
$C\langle0,1,1\rangle$.
Next, we show how to construct an equivalence.
Given a~Goursat bundle, an ef\/f\/icient method for constructing an equivalence map was worked out
in~\cite{Vassiliou2006b} and relies on the f\/iltration induced on the cotangent bundle.
Denote by~$\Xi^{(j)}$ and~$\Xi^{(j)}_{j-1}$ the annihilators of~$\ch\widehat{K}^{(j)}$ and
$\ch\widehat{K}^{(j)}_{j-1}$, respectively.
Then we obtain the f\/iltration
\begin{gather*}
\Xi^{(2)}\subset\Xi^{(2)}_1\subset\Xi^{(1)}
\end{gather*}
spanned as
\begin{gather*}
\{dq_1,dq_2,dq_3\}\subset\{dq_1,dq_2,dq_3,dx\}\subset\{dq_1,dq_2,dq_3,dx,dq_4,da_1\}.
\end{gather*}
The construction proceeds by building appropriate dif\/ferential operators and functions.
Because $m_k-m_{k-1}=1$, condition~2 of the def\/inition of Goursat manifold is vacuous.
Instead we f\/ix any f\/irst integral of $\ch\mathcal{D}^{(k-1)}$, denoted, $t$ and any section $Z$ of $\mathcal{D}$
such that $Zt=1$.
Then, def\/ine a~distribution $\Pi^k$ inductively as follows:
\begin{gather*}
\Pi^{\ell+1}=[Z,\Pi^\ell],\qquad \Pi^1=\ch\mathcal{D}^1_0,
\qquad
1\leq\ell\leq k-1.
\end{gather*}
There is a~function $\varphi^k$ which is a~f\/irst integral of~$\Pi^k$ such that $d\varphi^k\wedge dt\neq0$.
The function $\varphi^k$ is said to be a~{\it fundamental function of order}~$k$.
The space of fundamental functions of lower order can be constructed from the f\/iltration above by taking quotients.
Specif\/ically, as noted above, in this case the only fundamental functions of less than maximal order 3 are of order 2.
They are described by the quotient bundle
\begin{gather*}
\Xi^{(2)}_1/\Xi^{(2)}=\{dx\}\mod dq_1,dq_2,dq_3.
\end{gather*}
Without loss of generality we can take $\varphi^2=x$ to be a~fundamental function of order~2.
The construction of the equivalence map is now as follows.
Function $t$ is the ``independent variable" and successive dif\/ferentiation gives the higher order variables
\begin{gather*}
\begin{split}
& z^1=\varphi^2,
\qquad
z^1_1=Z\varphi^2,
\qquad
z^1_2=Z^2\varphi^2,
\\
& z^2=\varphi^k,
\qquad
z^2_1=Z\varphi^k,
\qquad
z^2_2=Z^2\varphi^k,
\qquad
z^2_3=Z^3\varphi^k.
\end{split}
\end{gather*}
We now implement this.
The f\/irst integrals of $\ch\widehat{K}^{(2)}$ are spanned by $q_1$, $q_2$, $q_3$ and any function of these can be chosen to
be $t$.
If we choose (say) $t=q_2$, then for $Z$ we choose
\begin{gather*}
Z=\frac{2}{q_2q_4}X,
\end{gather*}
where $X$ is the f\/irst vector f\/ield in the basis for $\widehat{K}$ above; for then $Zt=1$, as required.
We then construct the integrable distribution $\Pi^3$ as described above and discover that its f\/irst integrals are
spanned by
\begin{gather*}
q_2
\qquad
\text{and}
\qquad
\frac{q_1q_3}{1+q_1q_2}.
\end{gather*}
Hence the fundamental function of highest order is
\begin{gather*}
z^2=\varphi^k=\frac{q_1q_3}{1+q_1q_2}.
\end{gather*}
The data
\begin{gather*}
t=q_2,
\qquad
z^1=x,
\qquad
z^2=\frac{q_1q_3}{1+q_1q_2}
\end{gather*}
and dif\/ferentiation by $Z$ now constructs the local equivalence $\psi$ identifying $\widehat{K}=\widehat{H}_1$ and
$C\langle0,1,1\rangle$.
The local inverse $\psi^{-1}_1:\mathbb{R}\to\widehat{M}_1$ determines the integral submanifolds of $\widehat{H}_1$.

An exactly analogous calculation holds for $\widehat{H}_2$ and one arrives thereby at an explicit map
$\psi^{-1}_2:\mathbb{R}\to\widehat{M}_2$ representing the integral manifolds of $\widehat{H}_2$.
The explicit integral manifolds of $\Psi$ are a~{\it superposition} of those of $\widehat{H}_1$ and $\widehat{H}_2$:
\begin{gather*}
\mathbb{R}\times\mathbb{R}\to\boldsymbol{\sigma}\big(\psi^{-1}_1(\mathbb{R}),\psi^{-1}_2(\mathbb{R})\big).
\end{gather*}

In this way we obtain remarkably compact representations for wave maps into this metric:

\begin{theorem}[hyperbolic Weierstrass representation]
\label{Wrep}
For each collection of twice continuously differentiable real valued functions $f_1(s)$, $f_2(s)$, $g_1(t)$, $g_2(t)$
of parameters $s$, $t$, the functions
\begin{gather*}
x=f_1(s),\qquad y=g_1(t),
\\
e^{u/2}=
\frac{(tg_2(t)-1)\dot{f}_2(s)+(2tg_2(t)-1)f_2(s)^2}{g_2(t)(f_2(s)+s\dot{f}_2(s))},
\\
e^{v/2}=
\frac{(sf_2(s)-1)\dot{g}_2(t)+(2sf_2(s)-1)g_2(t)^2}{f_2(s)(g_2(t)+t\dot{g}_2(t))}
\end{gather*}
define harmonic maps
\begin{gather*}
\big(\mathbb{R}^{1,1}, dx dy\big)\to\left(N,\frac{du_1^2-du_2^2}{1+e^{-u_1}}\right)
\qquad
\text{by}
\qquad
u_1=u+v,
\qquad
u_2=u-v.
\end{gather*}
\end{theorem}

\section{Concluding remarks}\label{section5}

 We've seen that for a~large family of non-constant initial data, it is
possible to construct explicit integral representations for solutions of the Cauchy problem for wave maps into
a~certain non-constant curvature metric due to the fact that the corresponding Euler--Lagrange equation is Darboux
integrable and because the Cauchy data can be extended as a~f\/low by a~system of Lie type.
We have also constructed a~hyperbolic Weierstrass representation for such wave maps making use of the general theory
in~\cite{AFV,Vassiliou2006a,Vassiliou2006b}.
The fundamental ingredients throughout include the theory of systems of Lie type and the notion of a~Vessiot group
associated to any Darboux integrable exterior dif\/ferential system~\cite{AFV}.
We expressed the evolution of the Cauchy data as a~system of Lie type for the action of a~subgroup of the Vessiot group.

We may perhaps regret the occurence of {\it integral} representations in our solution of the Cauchy problem preferring
the elimination of all quadrature.
Alas, this is surely a~forlorn hope in a~dif\/f\/icult nonlinear problem, especially when it is recalled that even in
the general Cauchy problem for the $(1{+}1)$-linear wave equation quadrature cannot be eliminated, according to the
d'Alembert formula.
However, one can ask if there are Darboux integrable nonlinear sigma models with {\it solvable} Vessiot groups for harmonic
maps into nonzero curvature metrics.
This would make the application of the theory of systems of Lie type very useful indeed.
In fact there is at least one such sigma model~\cite{ClellandVassiliou}.

Interesting open problems include: what intrinsic properties of a~metric render the corresponding wave map system
Darboux integrable? Moreover, what do we learn about the geometry and topology of target manifolds from so vast
a~reduction in the Cauchy problem? This is especially intriging when it is recalled that the solution to the Cauchy
problem has been expressed as a~curve or ``evolution'' in a~f\/inite Lie group, arising by Lie's theory from
a~corresponding curve in its Lie algebra.

\appendix

\section{Appendix} \label{appendix}

\subsection{The superposition formula}\label{appendixA.1}

\looseness=-1
 For completeness, in this appendix we make
a~remark on the construction of the superposition formula $\boldsymbol{\sigma}:\widehat{M}_1\times\widehat{M}_2\to M$.
A general construction valid for any decomposable exterior dif\/ferential system was worked out in~\cite{AFV}.
Below we present results for the wave map system studied in this paper.
The hyperbolic structure of the wave map system in adapted coordinates is given by
\begin{gather*}
H=H_1\oplus H_2,
\end{gather*}
where
\begin{gather*}
H_1=\bigg\{\P x+a_2\P{a_1}+\frac{a_1(1+z_1z_2)}{2z_2z_4}\P{z_1}+\frac{z_2z_4}{2}\P{z_2}-\frac{a_1z_3}{2z_1z_2z_4}\P{z_3}+\frac{1}{2}(a_3z_4+a_1-z_4^2)\P{z_4},
\\
\phantom{H_1=\bigg\{}
\P{a_2},\, \P{a_3}\bigg\},
\\
H_2=\bigg\{\P y+b_2\P{b_1}+\frac{z_1z_3}{2}\P{z_1}+\frac{b_1(1+z_1z_2)}{2z_1z_3}\P{z_2}+\frac{1}{2}(b_3z_3-z_3^2+b_1)\P{z_3}-\frac{b_1z_4}{2z_1z_2z_3}\P{z_4},
\\
\phantom{H_2=\bigg\{}
\P{b_2},\, \P{b_3}\bigg\}.
\end{gather*}
Calculation shows that the inf\/initesimal symmetries of $H_1$ which are tangent to the level sets of all the f\/irst
integrals $x$, $y$, $\boldsymbol{a}$, $\boldsymbol{b}$, the tangential characteristic symmetries of~$H_1$, are spanned by
\begin{gather*}
\mathcal{E}_1=\bigg\{\frac{1}{2} (z_1z_3\P{z_1}-z_3^2\P{z_3} ),
\,
\frac{1+z_1z_2}{z_1z_3}\P{z_2}+\P{z_3}-\frac{z_4}{z_1z_2z_4}\P{z_4},
\,
-\frac{1}{2} (z_1\P{z_1}+z_2\P{z_2} )+z_3\P{z_3},
\\
\phantom{\mathcal{E}_1=\bigg\{}
{-}\frac{1}{4} (z_1\P{z_1}+z_2\P{z_2} )\bigg\}.
\end{gather*}
Similarly, the tangential characteristic symmetries of $H_2$ are spanned by
\begin{gather*}
\begin{split}
& \mathcal{E}_2=\bigg\{\frac{1+z_1z_2}{2z_2z_4}\P{z_1}-\frac{z_3}{z_1z_2z_4}\P{z_3}+\frac{1}{2}\P{z_4},
\\
& \phantom{\mathcal{E}_2=\bigg\{}
z_2z_4\P{z_2}-z_4^2\P{z_4},
\,
\frac{1}{2} (z_1\P{z_1}-z_2\P{z_2} )+z_4\P{z_4},
\,
-\frac{1}{2} (z_1\P{z_1}-z_2\P{z_2} )\bigg\}.
\end{split}
\end{gather*}
The structure of these Lie algebras are in {\it reciprocal relation}, namely, for $\mathcal{E}_1$ the nonzero Lie
brackets are
\begin{gather*}
[e_1,e_2]=e_3,\qquad [e_1,e_3]=e_1,\qquad [e_2,e_3]=-e_2,
\end{gather*}
while for $\mathcal{E}_2$ we have
\begin{gather*}
[f_1,f_2]=-f_3,\qquad [f_1,f_3]=-f_1,\qquad [f_2,f_3]=f_2
\qquad \mbox{and}
\qquad
[e_i,f_j]=0,
\quad
\forall
\,
i,j.
\end{gather*}
Hence, $\mathcal{E}_1$, $\mathcal{E}_2$ def\/ine the isomorphism class of the {\it Vessiot algebra} for this Darboux
integrable system, in the terminology of~\cite{AFV}.
Each $\mathcal{E}_i$ is isomorphic to $\mathfrak{sl}(2)\oplus\mathbb{R}$ and each frames a~neighbourhood of a~point of~$\mathbb{R}^4$.
By the converse of Lie's second fundamental theorem there is a~local Lie group $(G,m)$, where $m:G\times G\to G$ denotes group composition such that
$\mathcal{E}_i$ coincide with the inf\/initesimal left and right translations on~$G$.
From the expressions for either $\mathcal{E}_1$ or $\mathcal{E}_2$, by computing f\/lows or otherwise, we can compute
function $m$.
The components of $m$ coincide precisely with the f\/irst 4 components of the superposition map $\boldsymbol{\sigma}$,
where we interpret $\left(\{p_i\}_{i=1}^4,\{q_i\}_{i=1}^4\right)$ as local coordinates around $(e,e)\in G\times G$; $e$
being the identity in $G$.
\begin{remark} The {tangential characteristic symmetries}~\cite{Vassiliou2001,Vassiliou2000} described above do
not in general lead directly to the Vessiot algebra, which is needed for the construction of the superposition formula.
In general a~sequence of coframe adaptations described in~\cite{AFV} is required.
In terms of these the distribution $H_1\oplus H_2$ above can be constructed from components of the {\it fifth adapted
coframe} of a~so-called {\it Darboux pair}.
However, it turns out that for relatively low dimensional examples like the one in this paper, the full machinery
of~\cite{AFV} can sometimes be avoided and instead it is convenient to carry out a~direct computation of the tangential
characteristic symmetries of~$H_1$ and~$H_2$ as we did above.
A~simple example where this procedure is already not suf\/f\/icient is provided by a~Fermi--Pasta--Ulam equation
(see~\cite{Vassiliou2001}),
$u_{tt}=u_x^{-4}u_{xx}$.
\end{remark}

\subsection{Systems of Lie type}\label{appendixA2}

Let $\mu:G\times M\to M$ be a~left-action of a~Lie group $G$ on manifold $M$ and $\mathcal{G}\subset\mathfrak{X}(M)$
its Lie algebra of inf\/initesimal generators; $\mathfrak{X}(M)$ is the Lie algebra of all smooth vector f\/ields on~$M$.
Let $X_i$ be a~basis for~$\mathcal{G}$.
Then, a~{\it vector field of Lie type} or {\it Lie field} is a~curve in $\mathcal{G}$
\begin{gather*}
\mathcal{X}=\sum_i\alpha_i(t)X_i,
\end{gather*}
where $\alpha_i$ are smooth functions of parameter~$t$.
An {\it ODE of Lie type} is the dif\/ferential equation determined by a~Lie f\/ield
\begin{gather}
\label{LieEquation}
\frac{dx}{dt}=\mathcal{X}_{|_{x(t)}}.
\end{gather}

We pause to recall the map $\boldsymbol{\widehat{\mu}}:\mathfrak{g}\to\mathfrak{X}(M)$, def\/ined by
\begin{gather*}
\boldsymbol{\widehat{\mu}}(u)_{|_{x}}=\frac{d}{d\epsilon}\mu(\varphi_\epsilon(u),x)_{|_{\epsilon=0}},
\qquad
\forall\,
u\in\mathfrak{g},
\end{gather*}
where $\varphi_\epsilon(u)$ is the f\/low (in $G$) of $u\in\mathfrak{g}$ and $\mathfrak{g}$ is the matrix Lie algebra
associated to~$G$, viewed as a~matrix group.
The map $\boldsymbol{\widehat{\mu}}$ is the standard anti-homomorphism, induced by the left-action~$\mu$, from the
matrix Lie algebra $\mathfrak{g}$ to the Lie algebra of smooth vector f\/ields $\mathfrak{X}(M)$ on $M$.
Remark that for each $x\in M$, $\ker\boldsymbol{\widehat{\mu}}_{|_{x}}$ is equal to the subalgebra
$\mathcal{I}_x\subseteq\mathcal{G}$ of inf\/initesimal generators which vanish at $x$; that is,
$\ker{\boldsymbol{\widehat{\mu}}}_{|_{x}}$ is the {\it isotropy subalgebra} at $x$.
An action $\mu$ is {\it effective} on $M$ if the global isotopy group $\mathcal{I}(M)$ is trivial.
Recall that $\mathcal{I}(M)$ is the set of all $g\in G$ such that $\mu(g,x)=x$ for all $x\in M$.
It follows that if the action is ef\/fective then $\ker\boldsymbol{\widehat{\mu}}$ is trivial and the map
$\boldsymbol{\widehat{\mu}}$ is injective.
Henceforth we assume that the $G$-action is ef\/fective or at least almost ef\/fective; an action is {\it almost
effective}, if $\mathcal{I}(M)$ is discrete.
In this case, $\boldsymbol{\widehat{\mu}}$ is injective for those elements of $\mathfrak{g}$ whose f\/lows are close to
the the identity element in $G$.

We will also assume that the action is {\it transitive} so that for each $x,y\in M$ the corresponding isotropy
subalgebras are isomorphic, $\mathcal{I}_x\simeq\mathcal{I}_y$.
We therefore speak of {\it the } isotropy subalgebra or subgroup.
We now recall some useful notation from~\cite{DoubrovKomrakov}.

Let $f:M\to G$ be a~smooth map and let $\omega=\delta(f)=df\cdot f^{-1}$.
We know that $\omega$ is a~right-invariant Maurer--Cartan form and that it is valued in the Lie algebra $\mathfrak{g}$
of $G$.
\begin{definition}
Let $M$ be a~smooth manifold and $G$ a~Lie group.
\begin{enumerate}\itemsep=0pt
\item A dif\/ferential form $\omega\in\Omega(M,\mathfrak{g})$ is said to be a~Maurer--Cartan form if it satisf\/ies
\begin{gather*}
d\omega+\omega\wedge\omega=0.
\end{gather*}
\item A map $f:M\to G$ is said to be an {\it integral} of a~Maurer--Cartan form $\omega$ if and only if
$\delta(f)=\omega$.
\end{enumerate}
\end{definition}
\begin{theorem} On a~simply connected manifold every Maurer--Cartan form has an integral.
\end{theorem}

We denote the set of all $\mathfrak{g}$-valued Maurer--Cartan forms on $M$ by $\Omega(M,\mathfrak{g})$.
Specialising to the case $M=\mathbb{R}$ we have the following fundamental result.

\begin{theorem}
\label{fundamentalSol}
Let $\mu:G\times M\to M$ be an effective and transitive left-action of Lie group $G$ on manifold $M$.
Suppose $\omega=A(t)\,dt\in\Omega(\mathbb{R},\mathfrak{g})$ is smooth at $t=0$.
The unique solution of the initial value problem of Lie type
\begin{gather*}
\frac{dx}{dt}=\boldsymbol{\widehat{\mu}}(A(t))_{|_{x(t)}},
\qquad
x(0)=q\in M
\end{gather*}
is the function $x_q:\mathbb{R}\to M$ defined by
\begin{gather*}
x_q(t)=\mu(g(t),q),
\end{gather*}
where $g:\mathbb{R}\to G$, is the fundamental solution.
That is, $\delta(g)=\omega$ and $g(0)=e$.
\end{theorem}
\begin{remark} This theorem reduces the construction of the solution of a~system of Lie type to constructing the {\it
fundamental solution} $g:\mathbb{R}\to G$ satisfying
\begin{gather*}
\frac{dg}{dt}=A(t)g(t),
\qquad
g(0)=e.
\end{gather*}
This problem can still be very challenging.
Lie's approach is to simplify the problem by making use of any known solutions (Lie reduction).
See~\cite{DoubrovKomrakov} for a~proof of Theorem~\ref{fundamentalSol}.
\end{remark}

To explain this, continuing with transitive and ef\/fective left-action $\mu:G\times M\to M$, suppose that
$x_1(t),\ldots,x_k(t)$ are $k$ particular solutions of the Lie equation~\eqref{LieEquation}, satisfying
\begin{gather*}
x_j(t_0)=q_j,
\qquad
1\leq j\leq k,
\end{gather*}
for some collection of points $q_j\in M$.
Denote by $G_q$ the isotropy subgroup of $q\in M$,
\begin{gather*}
G_q=\{g\in G\, |\, g\cdot q=q\}.
\end{gather*}
Let
\begin{gather*}
H=G_{q_1}\cap G_{q_2}\cap\cdots\cap G_{q_k},
\end{gather*}
and let $\mathfrak{h}$ be the Lie algebra of $H$.
Furthermore, let a~curve $g_0:I\subseteq\mathbb{R}\to G$ be def\/ined by
\begin{gather}
\label{algebraicStep}
\mu(g_0(t),q_j):=g_0(t)\cdot q_j=x_j(t).
\end{gather}
Such a~curve of group elements is def\/ined up to a~multiplication on the right by a~curve in the the joint isotropy
subgroup of the initial conditions.
Constructing $g_0(t)$ involves the solution of algebraic equations with a~potentially large solution space with no
canonical choice of solution.
In general one aims to f\/ind a~solution $g_0(t)$ which passes through the identity at parameter value~$t_0$ where
initial conditions are to be posed.

Let $C^\infty(M,G)$ denote the set of smooth maps from $M$ to $G$.
As in~\cite{DoubrovKomrakov} def\/ine gauge transformations $\rho(h):\Omega(M,\mathfrak{g})\to\Omega(M,\mathfrak{g})$ by
\begin{gather*}
\rho(h)\omega=\text{Ad}(h)\omega+\delta(h),
\qquad
\forall\, h\in C^\infty(M,G).
\end{gather*}
\begin{theorem}[Lie reduction; see~\cite{Bryant, DoubrovKomrakov}]
\label{LieRed}
For any $\omega=A(t)\, dt\in\Omega(\mathbb{R},\mathfrak{g})$,
\begin{enumerate}\itemsep=0pt
\item[$1.$]
$
\omega_1=\rho\left(g_0(t)^{-1}\right)\omega\in\Omega(\mathbb{R},\mathfrak{h}).
$
\item[$2.$] If $g_1:I\subseteq\mathbb{R}\to H$ is an integral of $\omega_1$ then $g_0(t)g_1(t)$ is an integral of $\omega$.
\end{enumerate}
\end{theorem}

\subsubsection{Application of Lie reduction to the Riccati equation}
\label{LieRedRiccati}

The Cauchy problem for the integrable wave map system has been reduced to an initial value problem for Riccati
equation~\eqref{RiccatiIVP} together with one quadrature.
In this subsection we brief\/ly illustrate the use of {\it Lie reduction} in the resolution of the general Riccati
initial value problem.
It is used in Section~\ref{section3} of the paper.
The standard action of the special linear group $G={\rm SL}(2,\mathbb{R})$ on the real projective line induces
a~left-action on the real line
\begin{gather}
\label{Mobius}
{\lambda}_g(\xi):=\left[\begin{matrix}a&b\\ c&d\end{matrix}\right]\cdot\xi=\frac{a\xi+b}{c\xi+d},
\qquad
\xi\in\mathbb{R}
\end{gather}
the {\it linear fractional} or {\it M\"obius} transformations.

With basis
\begin{gather*}
u_{-1}=\left[\begin{matrix}0&1\\0&0\end{matrix}\right],
\qquad
u_0=\left[\begin{matrix}1&0\\0&-1\end{matrix}\right],
\qquad
u_1=\left[\begin{matrix}0&0\\1&0\end{matrix}\right]
\end{gather*}
for $\mathfrak{sl}(2,\mathbb{R})$ it is easy to see that the anti-homomorphism
$\widehat{\boldsymbol{\lambda}}:\mathfrak{sl}(2,\mathbb{R})\to\mathfrak{X}(\mathbb{R})$ can be expressed
\begin{gather*}
\widehat{\boldsymbol{\lambda}}(\alpha_0u_{-1}+\alpha_1u_0+\alpha_2u_1)=\big(\alpha_0+2\alpha_1x-\alpha_2x^2\big)\P x.
\end{gather*}
Hence
\begin{gather*}
A(t)=\alpha_0(t)u_{-1}+\alpha_1(t)u_0+\alpha_2(t)u_1=
\left[\begin{matrix}\alpha_1(t)&\alpha_0(t)\\ \alpha_2(t)&-\alpha_1(t)\end{matrix}\right]
\end{gather*}
satisf\/ies
\begin{gather*}
\widehat{\boldsymbol{\lambda}}\left(A(t)\right)=\big(\alpha_0(t)+2\alpha_1(t)x-\alpha_2(t)x^2\big)\P x.
\end{gather*}
That is, $A(t)\subset\mathfrak{sl}(2,\mathbb{R})$ is associated with the Lie equation
\begin{gather}
\label{riccati2}
\frac{dx}{dt}=\alpha_0(t)+2\alpha_1(t)x-\alpha_2(t)x^2.
\end{gather}
Now suppose that a~solution $x_0(t)$ of equation~\eqref{riccati2} is known and suppose that $x_0(0)=0$.
Solving equation~\eqref{algebraicStep} we obtain
\begin{gather*}
g_0(t)=\left[\begin{matrix}1&x_0(t)\\0&1\end{matrix}\right].
\end{gather*}
The isotropy subgroup of $0\in\mathbb{R}$ is
\begin{gather*}
H=G_0=\left[\begin{matrix}a&0\\ c&a^{-1}\end{matrix}\right],
\end{gather*}
from which the initial value problem for the fundamental solution $g_1(t)$ of the reduced equation of Lie type can be
deduced to be
\begin{gather}
\label{reducedFundamental}
\frac{dg_1}{dt}=B(t)g_1(t),
\qquad
g_1(0)=I_2,
\end{gather}
where
\begin{gather*}
B(t)dt=\delta(g_1(t))=\rho_{g_0(t)^{-1}}\left(\omega\right)=
\left[\begin{matrix}\alpha_1-x_0\alpha_2&0\\-\alpha_2&-(\alpha_1-x_0\alpha_2)\end{matrix}\right]dt.
\end{gather*}
Since $g_1(t)$ is a~curve in $H$, we have
\begin{gather*}
g_1(t)=\left[\begin{matrix}\gamma_1(t)&0\\\gamma_2(t)&\gamma_1(t)^{-1}\end{matrix}\right].
\end{gather*}
With $\beta(t)=\alpha_1-x_0\alpha_2$,~\eqref{reducedFundamental} has solution, expressed in terms of quadrature
\begin{gather*}
\gamma_1(t)=\exp\int_0^t~\beta(\tau)\,d\tau,
\qquad
\gamma_2(t)=-\gamma_1(t)^{-1}\int_0^t~\alpha_2(\tau)\gamma_1(\tau)^2\,d\tau.
\end{gather*}
Finally, the unique solution of the initial value problem $x(0)=q\in\mathbb{R}$ for the Riccati equation can now be
expressed in terms of quadrature:
\begin{gather*}
x(t)=g_0(t)g_1(t)\cdot q=
\left[\begin{matrix}\gamma_1(t)+x_0(t)\gamma_2(t)&x_0(t)/\gamma_1(t)\\
\gamma_2(t)&\gamma_1(t)^{-1}\end{matrix}\right]\cdot q.
\end{gather*}

That is,
\begin{gather*}
x(t)=\frac{\gamma_1(t)q}{\gamma_2(t)q+\gamma_1(t)^{-1}}+x_0(t).
\end{gather*}
\begin{theorem}
\label{quadReduction}
Let $x_0(t)$ be a~known solution of a~Riccati equation~\eqref{riccati2}.
Then by Lie reduction the construction of the fundamental solution is reducible to quadrature.
\end{theorem}

\subsection*{Acknowledgements}

I'm grateful to the three anonymous referees for their close reading of the manuscript and for making suggestions which
considerably improved the paper.
I would like to acknowledge, with my thanks, the early involvement of Jordane Math\'e for carefully working together
through the calculations in Section~\ref{section4} which formed a~portion of his internship from the Ecole normale sup\'erieure de
Cachan, France.
Much of the research for this paper was carried out while I was a~Visiting Fellow at the Mathematical Sciences
Institute of the Australian National University, Canberra.
The hospitality of the MSI is gratefully acknowledged.
In particular, I thank Mike Eastwood and the Dif\/ferential Geometry Group for stimulating discussions.

\pdfbookmark[1]{References}{ref}

\LastPageEnding

\end{document}